\documentclass[11pt,english]{article}
\usepackage[letterpaper]{geometry}
\usepackage{amsmath,amssymb,amsthm}
\usepackage{graphicx}
\usepackage{setspace}
\usepackage{babel}

\theoremstyle{plain}
\newtheorem{theorem}{Theorem}
\newtheorem*{theorem*}{Theorem}

\newtheorem*{corollary*}{Corollary}
\newtheorem{definition}{Definition}
\newtheorem*{remark*}{Remark}

\def\mathscr{\mathfrak}

\def\-{\backslash}

\makeatletter

\date{}

\makeatother

\begin{document}

\title{Bisecting horn angles}

\author{Sergiy Koshkin\\
 Department of Mathematics and Statistics\\
 University of Houston-Downtown\\
 One Main Street\\
 Houston, TX 77002\\
 e-mail: koshkins@uhd.edu}
\maketitle

\begin{abstract}\

A horn angle between a circle and its tangent is considered in Euclid's {\it Elements}, and Euclid remarks that it is smaller than any acute rectilinear angle. Already in antiquity, Proclus wondered whether it is possible to bisect horn angles. We will give a construction of a bisector which was within the means of ancient geometers since the time of Archimedes and Apollonius. We will also compare it to the conformal bisection method introduced in modern times. 
\bigskip

\bigskip


\end{abstract}

\section{History}\label{S1}

In book I of the {\it Elements}, Euclid defines angle as ``the inclination to one another of two lines which meet one another and do not lie in a straight line" \cite{HeEu}. ``Lines" in Euclid can be curves. So he allows curvilinear angles, but they are considered only once in the {\it Elements} \cite{Cab}. Namely, in proposition 16 of book III, Euclid characterizes a tangent to a circle as the line perpendicular to its diameter, and notes that ``the angle of the semicircle is greater, and the remaining angle less than any acute rectilineal angle". The ``remaining angle", the one between the circle and its tangent, became controversial already in antiquity. This is because Euclid's remark put it in tension with the definition of ratio, foundational to many parts of the {\it Elements} that anticipate what we now call real analysis. The definition, commonly attributed to Euclid's predecessor Eudoxus of Cnidus, reads: ``Magnitudes are said to have a ratio to one another which can when multiplied exceed one another". A horn angle, however, no matter how many times multiplied, does not exceed any acute rectilinear angle -- it is {\it infinitesimal}. Although this is presented as a definition, Euclid later wields it as an axiom, {\it assuming} that segments, areas, etc., do always have a ratio. This assumption also goes back to Eudoxus. Archimedes articulates it explicitly as axiom V in {\it On the Sphere and Cylinder}, and that all non-zero real numbers have ratios is now called the axiom of Archimedes. 

Euclid does not give a name to his ``remaining angle", but a 5-th century philosopher Proclus mentions one when commenting on the bisection of rectilinear angles. ``Bisection of any kind of angle", he writes, ``is not a matter for elementary treatise... Thus, it is difficult to say if it is possible to bisect so-called horn-like angle" \cite[I.9]{HeEu},  \cite{Proc}. A medieval mathematician Nemorarius called it ``the angle of contingence" (c. 1220), and Cardano ``the angle of contact" (1550). Its infinitesimal nature attracted the attention of many prominent later mathematicians, including Vieta, Galileo, Newton, Leibniz, and more recently, Hilbert and Klein \cite{Cab}. At the end of 19-th century, horn angles were placed into the general context of non-Archimedean analysis and geometry. In the 20-th century, the American mathematician Edward Kasner, and his students, studied the horn angles in the context of conformal geometry \cite{Kas12}-\cite{Pf14}. To readers interested in the history of horn angles, we recommend a series of 1996 posts by historian Jose Cabillon available online \cite{Cab}, which gives a comprehensive review, and includes an extensive bibliography of original sources. For a more recent work on  horn angles, see \cite{Asch, Bair, CiBu}.

Our present goal will be to do what Proclus found difficult to say could be done, to bisect the line-circle and circle-circle horn angles. First, we will do it in the spirit of the ancient Greeks, which will give us a new perspective on conic sections. Then we will briefly discuss Kasner's conformal approach. 

\section{Metric bisection}\label{S2}

\begin{figure}[htbp]
\vspace{-0.1in}
\begin{centering}
(a)\ \ \ \includegraphics[width=1.5in]{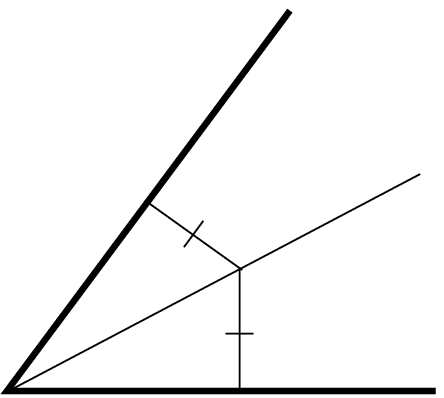} \hspace{1.0in} (b)\ \ \ \includegraphics[width=1in]{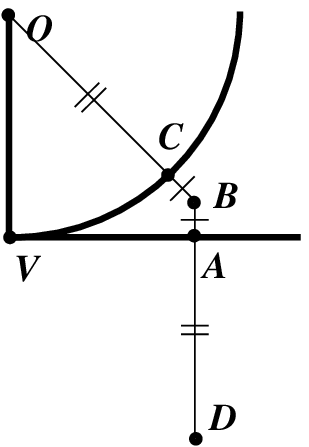}
\par\end{centering}
\vspace{-0.2in}
\hspace*{-0.1in}\caption{\label{MetrBis} (a) Rectilinear bisector; (b) Line-circle bisector}
\end{figure}

To bisect curvilinear angles, we first need to define what this means. Rectilinear angles are divided by a bisector into two congruent parts, but this will not work in general because the curved sides of a curvilinear angle are already not congruent. Another observation is that every point on a rectilinear bisector is equidistant from the sides of the angle, and this property specifies it uniquely, Figure \ref{MetrBis}(a). Following mathematical custom, we shall turn this property into a definition.
\begin{definition} A curve is called the metric bisector of a curvilinear angle if each point on it is equidistant from both sides of the angle.
\end{definition}
\noindent ``Metric" stands for the fact that our bisector is defined in terms of distances. Let us now consider Euclid's horn angle, between a line and a circle tangent to it. To find the set of points (locus, as the Greeks called it) equidistant to the line and the circle, one approach is to introduce coordinates and use the distance formula to derive its equation. But it is more illuminating to use Euclid style geometry instead, Figure \ref{MetrBis}(b).

\begin{figure}[htbp]
\vspace{-0.1in}
\begin{centering}
(a)\ \ \ \includegraphics[width=1in]{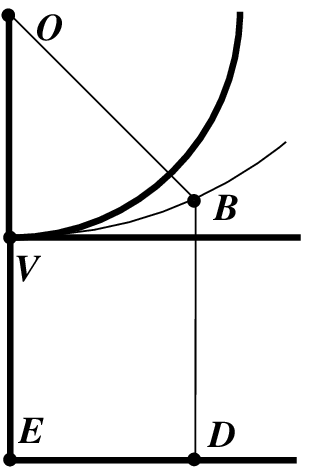} \hspace{1.0in} (b)\ \ \ \includegraphics[width=1.3in]{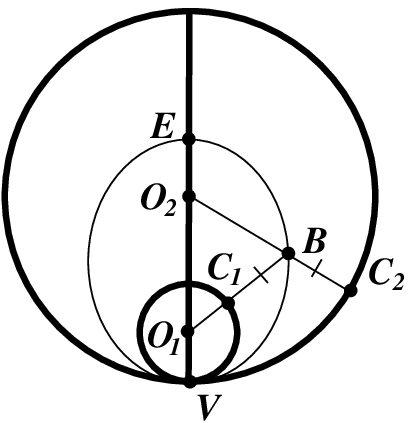}
\par\end{centering}
\vspace{-0.2in}
\hspace*{-0.1in}\caption{\label{CirCir} (a) Parabolic bisector with focus and directrix; (b) Ellipse as a circle-circle bisector}
\end{figure}

Note that the shortest segment $BC$ from a bisector point $B$ to the circle is perpendicular to the circle, and therefore lies on the line connecting $B$ to its center $O$. By definition, $AB=BC$, and $OC$ is of the same length $R$ for all points on the bisector, because it is the radius of the circle. So if we add a vertical segment $AD$ of length $R$ below $AB$ then also $OB=BD$. What is the point of adding $AD$, you ask? Well, points $D$ lie on a line parallel to $VA$, so this completes the picture geometrically, Figure \ref{CirCir}(a). Moreover, we now have characterized our locus as consisting of points equidistant from a point $O$, and a line $ED$. This is simpler than the original characterization, and is enough to tell us what the bisector is. The Greeks called such telling characterizations ``symptoms" \cite{Acer}.

As is known from precalculus texts, the locus equidistant from a point, called the {\it focus}, and a line, 
called the {\it directrix}, is a parabola. Nowadays, we often take this {\it focus-directrix property} as a geometric definition of parabolas. The distance between the focus and the vertex of a parabola is called its {\it focal distance}. Using these terms, we have proved the following theorem.
\begin{theorem} The metric bisector of a horn angle formed by a circle and a tangent is the parabola with focus at the center of the circle, focal distance equal to the circle's radius, and directrix parallel to the tangent.
\end{theorem}
The focus-directrix property was not the original definition of parabolas used by the ancient Greeks. The earliest known work where it appears is {\it On Burning Mirrors} by Diocles (c. 200 BC), who credits its discovery not to Archimedes, as many thought, but to his little known contemporary Dositheus, with whom Archimedes corresponded. Apparently, the Greeks did not think much of it, ``of this property, there is no trace in Apollonius's {\it Conica}" writes historian Fabio Acerbi \cite{Acer}, and the {\it Conica} is the paramount source for all things conic sections in antiquity. Even Diocles himself, writing about parabolic mirrors, buries it within an auxiliary construction.

Now let us turn to the circle-circle horn angles. There are two types of these angles, depending on whether the circles lie on the same side of their common tangent, or on opposite sides. We start with the first case, see Figure \ref{CirCir}(b).
By definition, $BC_1=BC_2$, and in this case, no auxiliary lines are necessary. Indeed, 
$BC_1=O_1B-O_1C_1=O_2C_2-O_2B=BC_2$. Therefore, $O_1B+O_2B=O_1C_1+O_2C_2$ is constant since the sum of the circle radii is the same for all points on the bisector. To determine the bisector, we need only to invoke another characterization known from precalculus texts: the locus of points, the sum of whose distances to two fixed points, called foci, remains constant, is an ellipse. This sum is equal to the major axis of the ellipse, because $O_1E+O_2E=O_1E+VO_1=VE$.

\begin{figure}[htbp]
\vspace{-0.1in}
\begin{centering}
(a)\ \ \ \includegraphics[width=1.2in]{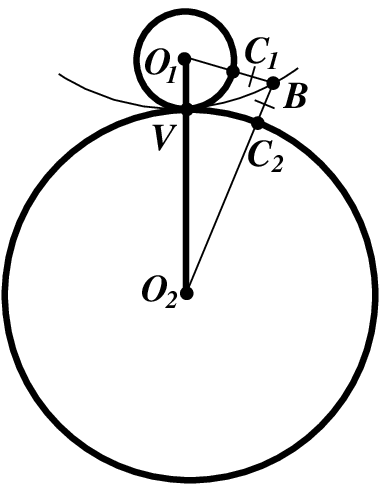} \hspace{0.5in} (b)\ \ \ \includegraphics[width=1.5in]{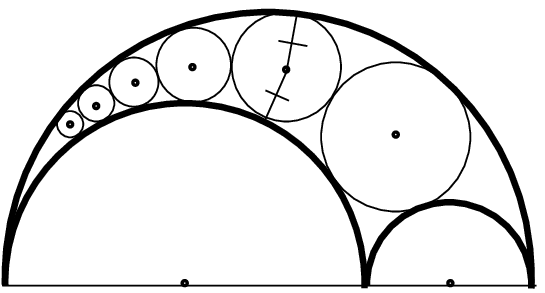}
\par\end{centering}
\vspace{-0.2in}
\hspace*{-0.1in}\caption{\label{Arb} (a) Hyperbola as a bisector; (b) Archimedes's arbelos and the Pappus chain}
\end{figure}

Unlike the focus-directrix property of a parabola, the focal property of an ellipse, and similar one of a hyperbola, with the difference of distances instead of the sum, was explicitly stated and proved in book 3 of Apollonius's {\it Conica}. This means that it was known earlier, in the first four books Apollonius streamlined and organized what was known before him, just as Euclid did in the {\it Elements}. We leave to the reader the case of circles lying on the opposite sides of a common tangent, see Figure \ref{Arb}(a). Summarizing, we have the following theorem.
\begin{theorem}\label{MetBis} The metric bisector of a horn angle formed by two circles lying on the same (opposite) side(s) of a common tangent is the ellipse (branch of the hyperbola) passing through the angle vertex, with foci at their centers, and the major axis equal to the sum (difference) of their radii.
\end{theorem}

Thus, all conic sections are neatly described as metric bisectors of line-circle and circle-circle horn angles. Why did ancient geometers not consider such constructions? In a way, they did. The elliptic bisector is closely related to the Pappus chain, a chain of successive tangent circles inscribed into Archimedes's arbelos, a shape formed by three semicircles sitting on the same segment, see Figure \ref{Arb}(b). Pappus, who lived in Alexandria in the 3rd century AD, left us the {\it Mathematical Collection}, one of the main surviving sources for ancient Greek mathematics, in which he proved that the centers of the circles in the chain lie on an ellipse \cite{Bank}. This is, of course, a direct consequence of Theorem \ref{MetBis}. But Pappus does not relate this result to bisecting a horn angle, nor does Proclus mention it in his discussion of them. Then again, Proclus says that ``it is not a matter for elementary treatise", i.e., not for the elementary straightedge and compass constructions, and the {\it Mathematical Collection} is a far more advanced treatise. 

\begin{figure}[htbp]
\vspace{-0.1in}
\begin{centering}
(a)\ \ \ \includegraphics[width=1.0in]{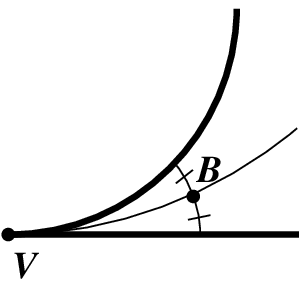} \hspace{0.5in} (b)\ \ \ \includegraphics[width=1.7in]{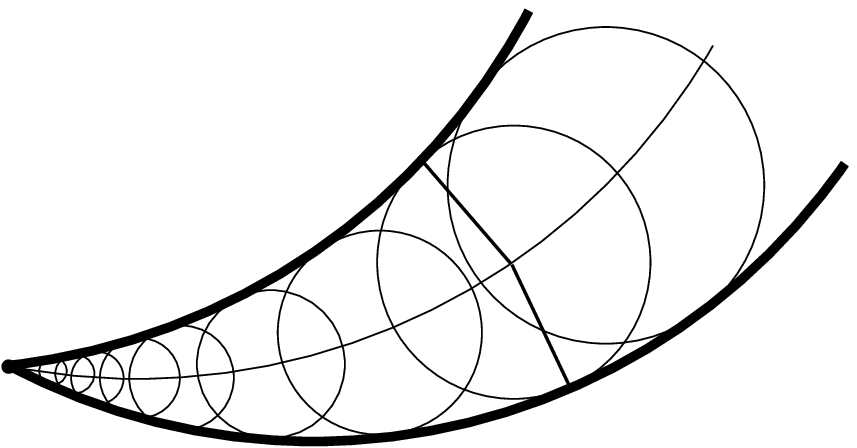}
\par\end{centering}
\vspace{-0.2in}
\hspace*{-0.1in}\caption{\label{Arc} (a) Bisector of circular arcs; (b) Bisector of wavefronts}
\end{figure}

Or perhaps, the reason was that there are different plausible definitions of bisecton that lead to different curves, and this is what Proclus meant by ``difficult to say if it is possible". For example, we could ask that every perpendicular through the bisector, or every arc of a circle centered at the vertex with the endpoints on the sides of the angle, be cut in half by the bisector, Figure \ref{Arc}(a). For rectilinear angles, these definitions give the same standard bisector, but for horn angles they do not. The reader may find it entertaining to play around with different definitions of bisection for lines and circles, and also to consider the inverse problem, of doubling the horn angle under them. 

The metric bisection has a nice physical interpretation, Figure \ref{Arc}(b). Imagine an object moving on the surface of water along a curve. A circular wave spreads from each point that it passes, and the outer edges of the disturbed water form, at any moment, a curvilinear angle, with the object at its vertex. Because the waves travel equal distances from the object to both edges, the path of the object is exactly the metric bisector of the curvilinear angle formed by them. Geometrically, the sides of the angle are the {\it envelopes} of a family of circles centered on the bisector, i.e., they are tangent to every circle in the family, a property that can be used to find the other side when one side and the bisector are known. In other words, this gives a way of metrically doubling a curvilinear angle.

\section{Conformal bisection}\label{S3}

At the end of the 17th century, Newton and Leibniz suggested using curvature for measuring horn angles \cite{Cab}, and at the end of the 19th century horn angles played a motivating role in the development of non-Archimedean geometry (Veronese, Hilbert) and analysis (du Bois-Reymond, Borel) \cite{Fish}. But there was not much research done on them until the American mathematician Edward Kasner suggested a new approach based on conformal geometry. His 1909-1912 results are summarized in a talk at the Fifth International Congress of Mathematicians \cite{Kas12}. As it applies to bisection, Kasner's idea was that a rectilinear bisector is exactly the line about which one side of the angle reflects into the other. The reflection about a line preserves both distances and angles, but one can not have that when reflecting curves -- different sides of a horn angle may not be congruent. However, if we settle for angles, such ``reflection" is possible to define when the sides are analytic curves, at least near the vertex. This was introduced by Hermann Schwarz at the end of 19th century, and plays an important role in the theory of analytic functions. However, it is not that well known, see \cite[\S 139]{Car}. 

\begin{figure}[htbp]
\vspace{-0.1in}
\begin{centering}
(a)\ \ \ \includegraphics[width=1.7in]{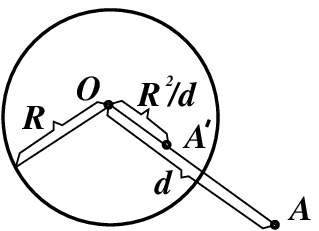} \hspace{1.0in} (b)\ \ \ \includegraphics[width=1.9in]{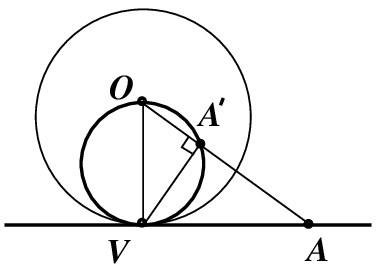}
\par\end{centering}
\vspace{-0.2in}
\hspace*{-0.1in}\caption{\label{Inv} (a) Circle inversion of a point; (b) Circle inversion of a tangent line}
\end{figure}

Fortunately, the Schwarz reflection about lines is the ordinary one, and about circles, it is the 
{\it circle inversion}, which is one of the M\"obius transformations known from complex analysis and projective geometry. For a circle of radius $R$, a point at a distance $d$ from the center is reflected to a point on the same ray through the center, but at a distance $R^2/d$, Figure \ref{Inv}(a). The outside of the circle is thereby reflected onto its inside (excluding the center, into which the ``infinity" is reflected), and the points on the circle are reflected into themselves. According to Pappus, already Apollonius considered the circle inversion in the lost work {\it On Plane Loci} (which is the Greeks' technical name for lines and circles). One of his results, reported in the {\it Mathematical Collection}, was that it transforms lines and circles into lines and circles \cite[Ch.3]{Rose}, a property common to all M\"obius transformations.

\begin{figure}[htbp]
\vspace{-0.1in}
\begin{centering}
(a)\ \ \ \includegraphics[width=1.2in]{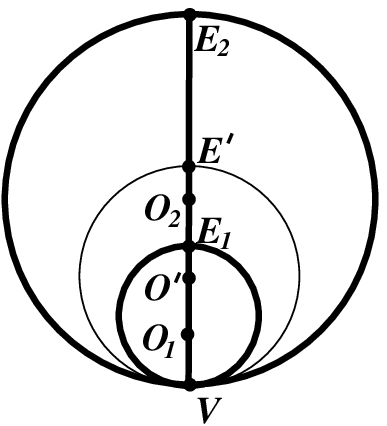} \hspace{1.0in} (b)\ \ \ \includegraphics[width=0.9in]{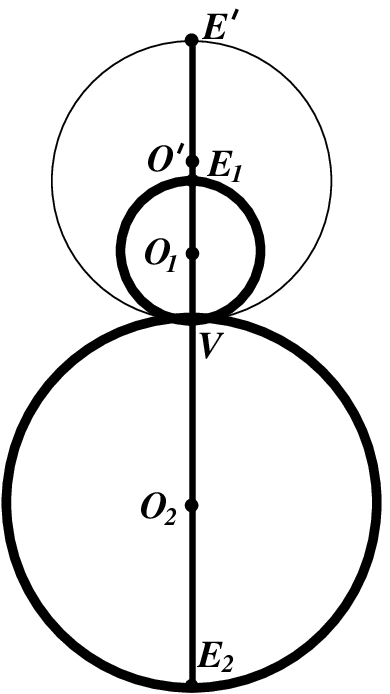}
\par\end{centering}
\vspace{-0.2in}
\hspace*{-0.1in}\caption{\label{ConfBis} Conformal circle-circle bisectors with (a) internal, (b) external tangency}
\end{figure}

Knowing the Schwarz reflection about a curve immediately solves the problem of doubling a horn angle, or any curvilinear angle -- we simply reflect one of its sides about the other. Figure \ref{Inv}(b) shows what happens when we reflect a tangent line to a circle about the circle. By definition of circle inversion, $OA\cdot OA'=OV^2$ or $OA:OV=OV:OA'$. This means that triangles $\triangle OV\!A$ and $\triangle OA'V$ have a common side ratio. Since they also have a common angle they are similar, and $\angle OA'V=\angle OV\!A$ is a right angle. But, as Euclid already knew, the locus of vertices of right 
angles subtending a line segment is the circle on it as a diameter. This follows from the converse to the inscribed angle theorem. In our case the segment is $OV$, the center of the circle is the midpoint of $OV$, and, therefore, it has half the radius of the reflecting circle. Now we are ready to turn to conformal bisections.
\begin{definition} A curve is called a conformal bisector of a curvilinear angle if one side of the angle is Schwarz  reflected about it into the other.
\end{definition}
\noindent 
As before, we will confine ourselves to the horn angles formed by lines and circles. A nice surprise is that not only are they themselves lines and circles, but they also can be constructed with straightedge and compass.
Reinterpreting Figure \ref{Inv}(b) as a conformal bisection of a line-circle horn angle $AV\!A'$, we see that the bisector is the circle tangent to both the line and the circle sides at the vertex, and is of twice the radius of the circle side. For a circle-circle angle with the circles on the same side of the tangent, the bisector has to be a circle tangent to both at the vertex. This is because the circle inversion reflects circles into circles, fixes points on the reflecting circle, and, being conformal, preserves the tangency, Figure  \ref{ConfBis}(a).
To find the bisector's radius, we note that the points $E_1$ and $E_2$, antipodal to the vertex, are reflections of each other about the bisector, so $O'E_1\cdot O'E_2=(O'E')^2$. Denoting by $R_1$ and $R_2$ the radii of the sides, and by $R'$ the radius of the bisector, we find:
$$
(2R_1-R')(2R_2-R')=R'^2,\text{ and therefore }R'=\frac{2R_1R_2}{R_1+R_2}=\frac{2}{\frac1{R_1}+\frac1{R_2}}.
$$ 
The right hand side is called the {\it harmonic mean} of $R_1$, $R_2$, this name goes back to the Pythagoreans. We can write this more suggestively if we use the curvatures $\kappa:=\frac1R$ of the circles instead of their radii. Then 
$\kappa':=\frac{\kappa_1+\kappa_2}2$, so the bisector's curvature is the arithmetic mean of the curvatures of its sides. The line-circle case is also covered if we take lines as circles of ``infinite radius", and $\kappa_2=0$. As before, we leave the case of circles on opposite sides of a tangent to the reader, and only state the end result: $\kappa':=\frac{\kappa_1-\kappa_2}2$, Figure \ref{ConfBis}(b). All three cases can be described uniformly in terms of {\it signed curvatures}, where we assign the curvature a positive sign if the circle's center is above the common tangent, and a negative sign, if it is below (``above" and ``below" can be chosen arbitrarily, but consistently for all circles). The general result is this.
\begin{theorem} The conformal bisector of a line-circle or a circle-circle horn angle is the circle (line) passing through the vertex, with the center on the line perpendicular to the common tangent, and the signed curvature equal to the arithmetic mean of the signed curvatures of the sides.
\end{theorem}
Similar reasoning allows us to conformally $n$-sect any line-circle or circle-circle horn angle. Since the radii of the $n$-sectors are rational in the radii of the sides, as in the case of the harmonic mean, they can be constructed with straightedge and compass \cite{Kas51}. In contrast, for rectilinear angles even trisection is already impossible with straightedge and compass.

As we saw, even for the simplest horn angles, the metric and the conformal bisectors are different curves. But do they match in some sense, at least near the vertex? In some sense, yes. For our original line-circle horn angle, Figure \ref{MetrBis}(b), let us place the origin at the vertex, and set the coordinate axes along the tangent and its perpendicular. Then the coordinate equation of the bisecting parabola is easily found to be $y=\frac{x^2}{4R}$. The calculus formula for the curvature of a graph gives $\kappa'(x):=\frac{\frac1{2R}}{\big(1+\left(\frac{x}{2R}\right)^2\big)^{3/2}}$, so $\kappa'(0)=\frac1{2R}=\frac{\kappa}2$, which matches the curvature of the conformally bisecting circle. With more calculations, the reader may verify that the same holds for the circle-circle bisectors. We believe that this is true in general, when both bisectors are defined, but proving it ``is not a matter for elementary treatise". The general problem of conformal bisection was studied by Kasner's student Pfeiffer in \cite{Pf14} and is subtle. A conformal bisector may not exist at all, but when it does, it is unique.

The results on conformal bisections suggest that Newton and Leibniz were right to use curvature for measuring horn angles. Taking the difference of curvatures at the vertex seems to give such a measure. But things are not so simple. A valid measure in conformal geometry has to be invariant under all conformal transformations. Transformations are called conformal if they preserve all (rectilinear) angles. But the difference of curvatures is not conformally invariant for general curves. Kasner showed that the simplest conformal invariant of a horn angle is the ratio $(\kappa_1-\kappa_2)^2:(\dot{\kappa_1}-\dot{\kappa_2})$, where $\dot{\kappa}$ denotes the derivative of the curvature with respect to arclength of the curve \cite{Kas37}. For lines and circles $\dot{\kappa}=0$, so the measure is determined by the difference of curvatures alone, as expected. More refined measures of curvilinear angles, that bring out the infinitesimal size of horn angles compared to rectilinear angles, are discussed in \cite{Asch,Bair}.
\bigskip

\noindent {\small {\bf Acknowledgements:} The author is grateful to the editor and the anonymous referee for the helpful comments that significantly improved the original version of the paper, and Linda Becerra for helping to prepare the final version.}

\end{document}